# MARTINGALE APPROXIMATIONS FOR SUMS OF STATIONARY PROCESSES

By Wei Biao Wu and Michael Woodroofe

*University of Chicago and University of Michigan*

Approximations to sums of stationary and ergodic sequences by martingales are investigated. Necessary and sufficient conditions for such sums to be asymptotically normal conditionally given the past up to time 0 are obtained. It is first shown that a martingale approximation is necessary for such normality and then that the sums are asymptotically normal if and only if the approximating martingales satisfy a Lindeberg–Feller condition. Using the explicit construction of the approximating martingales, a central limit theorem is derived for the sample means of linear processes. The conditions are not sufficient for the functional version of the central limit theorem. This is shown by an example, and a slightly stronger sufficient condition is given.

**1. Introduction.** The central limit problem for sums of stationary and ergodic processes has attracted continuing interest for over half a century, and two major lines of inquiry have developed. Under conditions of weak dependence such as strong mixing, blocking techniques have proved effective. Ibragimov (1962) provides an early account of this line. See Doukhan (1999) and Peligrad (1996) for more recent ones. An alternative approach, due to Gordin (1969), uses martingale approximation to establish asymptotic normality; see also Gordin and Lifsic (1978). Ho and Hsing (1997), Maxwell and Woodroofe (2000) and Wu and Woodroofe (2000) have followed this line recently. Here we come down on the side of martingale approximations by showing that if the partial sums of a stationary process are asymptotically normal in a suitable sense, then the martingale structure is present and the result could have obtained by using it. In addition, we sharpen the result of Maxwell and Woodroofe (2000), so that the necessary and sufficient conditions meet.

It is convenient to address the problem using the following notation. Let $(X_n)_{n\in\mathbb{Z}}$ be a stationary and ergodic Markov chain with values in the state space $\mathcal{X}$, and consider additive functionals

$$S_n = S_n(g) = \sum_{i=1}^{n} g(X_i), \tag{1}$$

where $g:\mathcal{X} \to \mathbb{R}$ is a measurable function for which $g(X_0)$ has mean 0 and finite variance. The partial sums of any stationary and ergodic process $(\xi_n)_{n\in\mathbb{Z}}$ may be written in this form by letting $X_n = (\ldots, \xi_{n-1}, \xi_n)$ and $g(X_n) = \xi_n$. Let $\pi$ denote the marginal distribution of $X_0$; suppose that there is a regular conditional distribution for $X_1$ given $X_0$, say $Q(x;B) = P(X_1 \in B | X_0 = x)$; and write $Qh(x) = E[h(X_1)|X_0 = x]$ a.e. $(\pi)$ for $h \in L^1(\pi)$.

Let $\sigma_n^2 = E(S_n^2)$ and suppose throughout the paper that

$$\sigma_n^2 \to \infty \tag{2}$$

as $n \to \infty$. This condition is needed to avoid degeneracy since otherwise there exists a stationary sequence $Y_n$ such that $g(X_n) = Y_n - Y_{n-1}$ [Theorem 18.2.2, Ibragimov and Linnik (1971)]. It will not be repeated in statements of lemmas and theorems. Consider a doubly indexed sequence $D_{nj}$ of random variables for which $D_{nj}$, $j = 1, 2, \ldots$, are martingale differences with respect to the filter $\mathcal{F}_j = \sigma(\ldots, X_{j-1}, X_j)$ for each $n$; and let $M_{nk} = D_{n1} + \cdots + D_{nk}$, so that $M_{nk}$, $k = 1, 2, \ldots$, is a martingale with respect to $\mathcal{F}_k$ for each $n$. The $D_{nj}$ or $M_{nk}$ is called a *martingale approximation* (to $S_n$) if

$$\max_{k \leq n} E[(S_k - M_{nk})^2] = o(\sigma_n^2). \tag{3}$$

A martingale approximation is called *stationary* if $D_{nj}$, $j = 1, 2, \ldots$, is a stationary sequence for each $n$, and *nontriangular* if $D_{nj} = D_j$ are independent of $n$. It is shown below that the existence of a stationary martingale approximation is equivalent to the existence of a nontriangular one. When (3) holds, asymptotic normality of $S_n/\sigma_n$ is equivalent to asymptotic normality of $M_{nn}/\sigma_n$, and this question may be addressed using the martingale central limit theorem [see, e.g., Billingsley (1995), pages 475–478].

It is shown in Section 2 that a simple growth condition on $E[E(S_n|X_0)^2]$ is necessary and sufficient for the existence of a martingale approximation. Then, in Section 3, it is shown that $S_n/\sigma_n$ is asymptotically standard normal, conditionally given $X_0$, iff the approximating martingales satisfy the conditions of the martingale central limit theorem. These conditions are not sufficient for the functional version of the central limit theorem. This is shown by example in Section 4, and a set of sufficient conditions is developed there.

Dedecker and Merlevede (2002) have used blocking techniques to obtain necessary and sufficient conditions for conditional asymptotic normality without assuming that the process is strongly mixing, or even ergodic. One of their conditions is closely related to (4), but their conditions do not include the existence of a martingale approximation and their uniform integrability condition for $S_n^2/n$ looks quite different from our Lindeberg–Feller conditions, (11) and (12). Using the explicit construction of martingales, we are able to obtain novel asymptotic theory for the sample means of linear processes, important and widely used stationary processes.

**2. Martingale approximations.** Below, $\|\cdot\|$ denotes the norm in an $L^2$ space, which may vary from one use to the next. For example, $\|\cdot\|$ denotes the norm in $L^2(P)$ in (4), and the norm in $L^2(\pi)$ in (5).

LEMMA 1. *If*

$$\|E(S_n|X_0)\| = o(\sigma_n), \tag{4}$$

*as $n \to \infty$, then there is a slowly varying function $\ell$ for which $\sigma_n^2 = n\ell(n)$.*

PROOF. If relation (4) holds, then $|E[S_n(S_{n+m} - S_n)]| = |E[S_n E(S_{n+m} - S_n | X_n)]| \leq \|S_n\| \times \|E(S_m|X_0)\| \leq \varepsilon_m \sigma_m \sigma_n$, where $\varepsilon_m \to 0$ as $m \to \infty$. The lemma then follows from Ibragimov and Linnik [(1971), Theorem 18.2.3 and the Remark on page 330], after correcting for obvious typographical errors. □

Relation (4) is crucial in what follows. Since $\|E[g(X_k)|X_0]\| = \|Q^k g\|$, it is implied by the condition, $\sum_{k=1}^n \|E[g(X_k)|X_0]\| = \sum_{k=1}^n \|Q^k g\| = o(\sigma_n)$, on the individual summands; but (4) is weaker and not unintuitive.

Recall that the equation $h = Qh + g$ is called *Poisson's equation*. Below, we will call a sequence $h_n \in L^2(\pi)$ an *approximate solution to Poisson's equation (for $g$)* if

$$\|h_n\| + n\|(I - Q)h_n - g\| = o(\sigma_n) \tag{5}$$

as $n \to \infty$. Also, if $a_n$ and $b_n$ are positive sequences, then $a_n \sim b_n$ iff $\lim_{n\to\infty} a_n/b_n = 1$.

THEOREM 1. *The following are equivalent*:

(i) *Relation* (4) *holds.*
(ii) *There is an approximate solution to Poisson's equation* (5).
(iii) *There is a stationary martingale approximation* (3).
(iv) *There is a nontriangular martingale approximation.*

In this case $E(D_{n1}^2) \sim \ell(n)$ for any stationary martingale approximation; and there is a stationary martingale approximation for which $\max_{k \leq n} \|S_k - M_{nk}\| \leq 3 \max_{k \leq n} \|E(S_k|X_0)\|$.

PROOF. It will be shown first that (i) $\Rightarrow$ (ii) $\Rightarrow$ (iii) $\Rightarrow$ (i) and then that (iii) $\Rightarrow$ (iv) $\Rightarrow$ (i). The remainder of the proof is placed between the two equivalences.

(i) $\Rightarrow$ (ii). If (4) holds, let $h_n^o = g + Qg + \cdots + Q^{n-1}g$. Then $h_n^o(x) = E(S_n|X_1 = x)$ and $Qh_n^o(x) = E(S_n|X_0 = x)$ for a.e. $x$. Clearly, $h_n^o = g + Qh_n^o - Q^n g$, and $\|h_n^o - Qh_n^o\| \leq 2\|g\|$. So, $\|h_n^o\| \leq \|E(S_n|X_0)\| + 2\|g\| = o(\sigma_n)$, by (4). Next, let

$$h_n = \frac{h_1^o + \cdots + h_n^o}{n}.$$

Then $\|h_n\| \leq \max_{k \leq n} \|h_k^o\| = o(\sigma_n)$, and $h_n = g + Qh_n - Qh_n^o/n$. So, $n\|(I - Q)h_n - g\| \leq \|Qh_n^o\| = o(\sigma_n)$, establishing (5).

(ii) $\Rightarrow$ (iii). If (5) holds, let $f_n = g - (I - Q)h_n$,

(6) $$D_{nk} = h_n(X_k) - Qh_n(X_{k-1}),$$

and $M_{nk} = D_{n1} + \cdots + D_{nk}$ for $k = 1, 2, \ldots$. Then $D_{n1}, D_{n2}, \ldots$ are stationary martingale differences for each $n$. Next, writing $g(X_k) = h_n(X_k) - Qh_n(X_k) + f_n(X_k)$ in (1) and rearranging terms then leads to $S_k = M_{nk} + S_k(f_n) + Qh_n(X_0) - Qh_n(X_k)$. So,

$$\max_{k \leq n} \|S_k - M_{nk}\| \leq n\|f_n\| + 2\|Qh_n\| = o(\sigma_n),$$

and (3) holds.

(iii) $\Rightarrow$ (i). If (3) holds, then $\|E(S_n|X_0)\| = \|E(S_n - M_{nn}|X_0)\| \leq \|S_n - M_{nn}\| = o(\sigma_n)$. This establishes the equivalence of (i)–(iii).

For any stationary martingale approximation in (3), $nE(D_{n1}^2) = E(M_{nn}^2) = E(S_n^2) + o(\sigma_n^2) \sim \sigma_n^2$, so that $E(D_{n1}^2) \sim \ell(n)$; and for the stationary martingale approximation constructed in the proof of (i) $\Rightarrow$ (iii), $\max_{k \leq n} \|S_k - M_{nk}\| \leq n\|f_n\| + 2\|Qh_n\| \leq 3 \max_{k \leq n} \|Qh_k^o\| \leq 3 \max_{k \leq n} \|E(S_n|X_0)\|$.

(iii) $\Rightarrow$ (iv) $\Rightarrow$ (i). If there is a stationary martingale approximation, then (4) holds and there is a stationary martingale approximation of the form (6), say $M_{nk} = D_{n1} + \cdots + D_{nk}$. Then $\|M_{nk} - M_{mk}\| \leq \|S_k - M_{nk}\| + \|S_k - M_{mk}\|$, and $m\|D_{n1} - D_{m1}\|^2 = \|M_{nm} - M_{mm}\|^2 \leq 2\|S_m - M_{mm}\|^2 + 2\|S_m - M_{nm}\|^2$. Let $D_k = D_{kk}$ and $M_n = D_1 + \cdots + D_n$. Then $M_1, M_2, \ldots$ is a martingale, $\|S_n - M_n\| \leq \|S_n - M_{nn}\| + \|M_{nn} - M_n\|$, and $\|S_n - M_{nn}\| = o(\sigma_n)$, by assumption. Here $\|M_{nn} - M_n\|^2 = \sum_{k=1}^n \|D_{nk} - D_{kk}\|^2 = \sum_{k=1}^n \|D_{n1} - D_{k1}\|^2$. So,

$$\|M_{nn} - M_n\|^2 \leq \sum_{k=1}^n \frac{2\|S_k - M_{kk}\|^2}{k} + \sum_{k=1}^n \frac{2\|S_k - M_{nk}\|^2}{k} = I_n + II_n,$$

say. Karamata's theorem [see, e.g., Theorem 0.6 in Resnick (1987)] implies that, for $\alpha > -1$, $\sum_{i=1}^{n} i^\alpha \ell(i) \sim n^{1+\alpha}\ell(n)/(1+\alpha)$. Hence $I_n = o[\sum_{k=1}^{n} \ell(k)] = o[n\ell(n)] = o(\sigma_n^2)$. For the second term, notice that $\|M_{nk}\|^2 = k\|D_{n1}\|^2$ and $\|D_{n1}\|^2 \sim \ell(n)$. Then for some positive $C$ and any positive $\varepsilon < 1/2$,

$$II_n \leq 4 \sum_{k \leq n\varepsilon} \frac{\|S_k\|^2 + \|M_{nk}\|^2}{k} + 2 \sum_{n\varepsilon < k \leq n} \frac{\|S_k - M_{nk}\|^2}{k}$$

$$\leq C \sum_{k \leq n\varepsilon} [\ell(k) + \ell(n)] + \frac{2}{\varepsilon} \max_{k \leq n} \|S_k - M_{nk}\|^2,$$

which by Karamata's theorem implies that $\limsup_{n \to \infty} II_n/\sigma_n^2 \leq 2C\varepsilon$ and, therefore, $\limsup_{n \to \infty} II_n/\sigma_n^2 = 0$. Conversely, if there is a nontriangular martingale approximation, then $\|E(S_n|X_0)\| = \|E(S_n - M_n|X_0)\| \leq \|S_n - M_n\| = o(\sigma_n)$, as above. □

As it is clear from Theorem 1, martingale approximations are not unique. Any two are asymptotically equivalent, however, in the following sense: If (3) holds, and if $M'_{nk} = D'_{n1} + \cdots + D'_{nk}$ is a second martingale approximation, then

(7) $\quad E\left[\max_{k \leq n}(M'_{nk} - M_{nk})^2\right] \leq 4\|M'_{nn} - M_{nn}\|^2 = 4\sum_{k=1}^{n} \|D'_{nk} - D_{nk}\|^2,$

using Doob's [(1953), page 317] inequality, and $\|M'_{nn} - M_{nn}\| \leq \|S_n - M'_{nn}\| + \|S_n - M_{nn}\| = o(\sigma_n)$.

If $\ell(n) \to \infty$ in Lemma 1, then it is impossible to have a martingale approximation that is both nontriangular and stationary, but if $\sigma_n^2 \sim \sigma^2 n$, then it is. Maxwell and Woodroofe (2000) show that if $\sum_{n=1}^{\infty} n^{-3/2}\|E(S_n|X_0)\| < \infty$, then there is a martingale $M_1, M_2, \ldots$ with stationary increments for which $\|S_n - M_n\|^2 = o(n)$. A simplified proof of a special case of this result is provided in Lemma 5, along with an explicit bound on $\|S_n - M_n\|$.

The proof of Theorem 1 contains the explicit construction of $D_{nk} = h_n(X_k) - Qh_n(X_{k-1})$ in terms of any approximate solution $h_n$ to Poisson's equation and also an explicit construction of $h_n$. An alternative approximate solution to Poisson's equation is provided next.

COROLLARY 1. *If* (4) *holds, then* (5) *holds with* $h_n = f_{1/n}$, *where*

$$f_\varepsilon(x) = \sum_{j=1}^{\infty}(1+\varepsilon)^{-j}Q^{j-1}g$$

*for* $0 < \varepsilon < 1$.

PROOF. From the definition, it is clear that $(1+\varepsilon)f_\varepsilon = g + Qf_\varepsilon$ and $(I-Q)h_n = g - h_n/n$. So, the corollary would follow from $\|h_n\| = o(\sigma_n)$. To see this, first observe that $f_\varepsilon = \varepsilon \sum_{k=1}^\infty (1+\varepsilon)^{-k-1} h_k^o$, by partial summation, where $h_k^o(x) = E(S_k|X_1 = x)$, as above. Let $V(s) = \sum_{k=1}^\infty \sigma_k s^{k+1}$. Then $\|h_n\| = o[V(n/(n+1))]/n$ by (4), and $V(s) \sim \frac{1}{2}\sqrt{\pi}(1-s)^{-3/2}\ell^{1/2}(1/(1-s))$ as $s \uparrow 1$ by Tauberian's theorem [see, e.g., Feller (1971), page 445]. □

For some examples, let $\ldots \eta_{-1}, \eta_0, \eta_1, \ldots$ be a stationary sequence of martingale differences with finite variance; and let $\ldots \theta_{-1}, \theta_0, \theta_1, \ldots$ be a sequence of i.i.d. random elements that is independent of $\ldots \eta_{-1}, \eta_0, \eta_1, \ldots$. Then $X_k = [(\ldots \theta_{k-1}, \theta_k), (\ldots, \eta_{k-1}, \eta_k)]$ is a stationary Markov chain with values in $\mathcal{X} = \Theta^{\mathbb{N}} \times \mathbb{R}^{\mathbb{N}}$, where $\Theta$ is the range of the $\theta_k$ and $\mathbb{N}$ is the nonnegative integers. Let $a_j : \mathcal{X} \to R$ be measurable functions for which

$$\sum_{j=0}^\infty E[a_j(X_0)^2 \eta_1^2] < \infty.$$

Then

(8) $$\xi_k = \sum_{j=0}^\infty a_j(X_{k-j-1})\eta_{k-j}$$

converges w.p.1 for each $k$ and is of the form $g(X_k)$. Processes of the form (8) include linear processes with constant $a_j$ and $\theta_k \equiv 0$, and are called *quasi-linear processes* below. They also include many nonlinear time series models, like autoregressive processes with random coefficients. Writing $\xi_k = \sum_{j \leq k} a_{k-j}(X_{j-1})\eta_j$ and letting $b_n = a_0 + \cdots + a_n$, it is easily seen that

$$E(S_n|X_0) = \sum_{j \leq 0}[b_{n-j}(X_{j-1}) - b_{-j}(X_{j-1})]\eta_j,$$

$$S_n - E(S_n|X_0) = \sum_{j=1}^n b_{n-j}(X_{j-1})\eta_j.$$

So, $\sigma_n^2 = \sigma_{n,1}^2 + \sigma_{n,2}^2$, with

$$\sigma_{n,1}^2 = \|E(S_n|X_0)\|^2 = \sum_{j=0}^\infty E\{[b_{j+n}(X_0) - b_j(X_0)]^2 \eta_1^2\},$$

$$\sigma_{n,2}^2 = \|S_n - E(S_n|X_0)\|^2 = \sum_{j=1}^{n-1} E[b_j(X_0)^2 \eta_1^2],$$

and (4) is equivalent to $\sigma_{n,1}^2 = o(\sigma_{n,2}^2)$. In this case, by (6) in the proof of Theorem 1, $D_{nk} = \bar{b}_n(X_{k-1})\eta_k$, where $\bar{b}_n = (b_0 + \cdots + b_{n-1})/n$, by some routine calculations, and $E[\bar{b}_n(X_0)^2 \eta_1^2]$ must be slowly varying. Observe that if $b_j(X_0)$, $j \leq 0$, are independent of $\eta_1$, then $E[\bar{b}_n(X_0)^2 \eta_1^2] = E[\bar{b}_n(X_0)^2] E(\eta_1^2)$ and that $\sigma_{n,1}^2$ and $\sigma_{n,2}^2$ simplify similarly.

EXAMPLE 1 (Linear processes). Suppose that $a_n$ are constants and (without loss of generality) that $E(\eta_k^2) = 1$. Then $b_n$ are also constants, and $\sigma_{n,1}^2 = o(\sigma_{n,2}^2)$ iff

$$(9) \qquad \sum_{j=0}^{\infty} (b_{j+n} - b_j)^2 = o\left[\sum_{k=1}^{n-1} b_k^2\right].$$

If $a_n$ are absolutely summable and $b := \sum_{n=0}^{\infty} a_n \neq 0$, then $\sigma_{n,2}^2 \sim b^2 n$ and there is a $C$ for which $\sigma_{n,1}^2 \leq C \sum_{i=1}^{n} \sum_{j=i}^{\infty} |a_j| = o(n)$, so that (9) holds. Relation (9) also holds if $b \neq 0$ and $b_n = b + O(1/n)$. If $a_0 = 0$ and $a_n = 1/n$ for $n \geq 1$, then $b_n \sim \log(n)$, and $\sigma_{n,2}^2 \sim n \log^2(n)$. In this case $\sigma_{n,1}^2 = O(n) = o(\sigma_{n,2}^2)$, so that (9) holds. To see this, observe that, for $j \geq 3$, $1/(j+1) \leq \int_j^{j+1} u^{-1} du$ and $[\log(j+n) - \log(j)]^2 \leq \int_{j-1}^{j} [\log(u+n) - \log(u)]^2 du$, so that

$$\sum_{j=3}^{\infty} (b_{j+n} - b_j)^2 \leq \sum_{j=3}^{\infty} \left(\int_j^{j+n} \frac{1}{u} du\right)^2$$

$$= \sum_{j=3}^{\infty} \log^2 \frac{j+n}{j}$$

$$\leq \int_2^{\infty} \log^2 \frac{u+n}{u} du = O(n).$$

Similarly, if $a_0 = 0$, $a_1 = 1/\log(2)$ and $a_n = 1/\log(n+1) - 1/\log(n)$ for $n \geq 2$, then $\sigma_{n,2}^2 \sim n/\log^2(n)$ and $\sigma_{n,1}^2 = O[n/\log^3(n)] = o(\sigma_{n,2}^2)$, so that (4) holds. On the other hand, if $a_n = n^{-\beta}$, where $1/2 < \beta < 1$, then there are positive constants $c_{1,\beta}$ and $c_{2,\beta}$ for which $\sigma_{n,i}^2 \sim c_{i,\beta} n^{3-2\beta}$ as $n \to \infty$ for $i = 1, 2$, so that (4) fails.

**3. Asymptotic normality.** The main result of this section is that $S_n^* := S_n/\sigma_n$ is asymptotically standard normal given $X_0$, as described below, iff there is a martingale approximation, (3), and the $D_{nk}$ satisfy the conditions of the martingale central limit theorem, (11) and (12). In more detail, let $P^x$ and $E^x$ denote the regular conditional probability and conditional expectation for $\mathcal{F}_\infty$ given $X_0 = x$; and let $F_n$ denote the conditional distribution function

$$F_n(x; z) = P^x(S_n^* \leq z).$$

Further let $\Phi$ denote the standard normal distribution function; and let $\Delta$ denote the Levy distance between two distribution functions. Then by asymptotic normality given $X_0$, we mean

$$\lim_{n \to \infty} \int_{\mathcal{X}} \Delta[\Phi, F_n(x; \cdot)] \pi\{dx\} = 0. \tag{10}$$

Clearly, (10) implies that $S_n^*$ is asymptotically standard normal, but (10) is stronger in general; it implies that $S_n^*$ is asymptotically standard normal for related models in which $X_0$ has any distribution that is absolutely continuous with respect to the stationary distribution. Such a property is needed, for example, if asymptotic normality is used to set approximate error bounds for Markov chain Monte Carlo experiments. See, for example, Tierney (1994). Under conditions of weak dependence, (10) can be deduced from (unconditional) asymptotic normality of $S_n^*$. See Proposition 1 for the details and the continuation of Example 1 for a case in which $S_n^*$ is (unconditionally) normal, but (10) fails.

LEMMA 2. *If* (10) *holds, then* (4) *holds; that is,* $\|E(S_n|X_0)\| = o(\sigma_n)$.

PROOF. The proof follows Maxwell and Woodroofe (2000), who considered the special case $\sigma_n^2 \sim cn$; it is included because the lemma is crucial to what follows. Let $\Rightarrow$ denote convergence in distribution. Notice that if $Z_m \Rightarrow \Phi$, then $\liminf_{m \to \infty} \text{var}(Z_m) \geq 1$, where $\text{var}(Z_m) = E(Z_m^2) - [E(Z_m)]^2$. To see this, for $J > 0$, let $T_{m,J} = \min[\max(Z_m, -J), J]$. Then $\lim_{m \to \infty} \text{var}(T_{m,J}) = \int_{\mathbb{R}} \{\min[\max(u, -J), J]\}^2 \, d\Phi(u) \underset{J \to \infty}{\to} 1$. By Corollary 4.3.2 in Chow and Teicher (1978), $\text{var}(Z_m) \geq \text{var}(T_{m,J})$. So $\liminf_{m \to \infty} \text{var}(Z_m) \geq 1$.

Assume otherwise that there is a $\delta > 0$ such that $\|\mathbb{E}(S_{n'}^*|X_0)\| > \delta$ along a subsequence $\{n'\}$. By (10), there exists a further subsequence $\{n''\} \subset \{n'\}$ such that $\Delta[\Phi, F_{n''}(x; \cdot)] \to 0$ for almost all $x(\pi)$. Let $\tau_n^2(x) = \text{var}(S_n^*|X_0 = x)$. By the result in the previous paragraph, $\liminf_{n'' \to \infty} \tau_{n''}^2(x) \geq 1$ for almost all $x(\pi)$. Thus $1 \leq \liminf_{n'' \to \infty} \int_{\mathcal{X}} \tau_{n''}^2(x) \pi(dx)$ by Fatou's lemma. On the other hand, the integral in the previous inequality equals $\|S_{n''}^*\|^2 - \|\mathbb{E}(S_{n''}^*|X_0)\|^2 \leq 1 - \delta^2$, which is a contradiction. $\square$

LEMMA 3. *Suppose there is a martingale approximation* $\{D_{nk}\}$ *for which*

$$\frac{1}{\sigma_n^2} \sum_{k=1}^{n} E(D_{nk}^2 | \mathcal{F}_{k-1}) \xrightarrow{p} 1 \tag{11}$$

*and*

$$\frac{1}{\sigma_n^2} \sum_{k=1}^{n} E(D_{nk}^2 \mathbf{1}_{\{|D_{nk}| \geq \varepsilon \sigma_n\}} | \mathcal{F}_{k-1}) \xrightarrow{p} 0 \tag{12}$$

hold for each $\varepsilon > 0$. Then for any martingale approximation $\{D'_{nk}\}$ (say), (11) and (12) are satisfied. In addition,

$$(13) \qquad \sup_{0 < t \leq 1} \left| \frac{1}{\sigma_n^2} \sum_{k \leq nt} E(D'^2_{nk}|\mathcal{F}_{k-1}) - t \right| \xrightarrow{p} 0.$$

PROOF. Observe that $E|E(D'^2_{nk}|\mathcal{F}_{k-1}) - E(D^2_{nk}|\mathcal{F}_{k-1})| \leq E|D'^2_{nk} - D^2_{nk}|$ and

$$E(D'^2_{nk} \mathbf{1}_{\{|D'_{nk}| \geq 2\varepsilon\sigma_n\}}|\mathcal{F}_{k-1})$$
$$\leq 2E(D^2_{nk} \mathbf{1}_{\{|D_{nk}| \geq \varepsilon\sigma_n\}}|\mathcal{F}_{k-1}) + 2E(|D'^2_{nk} - D^2_{nk}||\mathcal{F}_{k-1}).$$

So, if $D_{nk}$ satisfies (11) and (12), then so do $D'_{nk}$, since

$$E\left(\sum_{k=1}^{n} |D'^2_{nk} - D^2_{nk}|\right) \leq \sqrt{\sum_{k=1}^{n} \|D'_{nk} + D_{nk}\|^2} \times \sqrt{\sum_{k=1}^{n} \|D'_{nk} - D_{nk}\|^2} = o(\sigma_n^2),$$

as in (7). To establish (13), let $m = \lfloor nt \rfloor$, where $\lfloor x \rfloor$ is the greatest integer that does not exceed $x$; let $M'_{nk} = D'_{n1} + \cdots + D'_{nk}$. Observe that $\sigma_m^2/\sigma_n^2 \to t$ as $n \to \infty$, (11) implies

$$\frac{1}{\sigma_n^2} \sum_{k=1}^{m} E(D'^2_{mk}|\mathcal{F}_{k-1}) \xrightarrow{p} t.$$

Since $\|M'_{nm} - M'_{mm}\| \leq \|M'_{nm} - S_m\| + \|S_m - M'_{mm}\| = o(\sigma_n)$,

$$E\left(\sum_{k=1}^{m} |D'^2_{nk} - D'^2_{mk}|\right) \leq \sqrt{\sum_{k=1}^{m} \|D'_{nk} + D'_{mk}\|^2} \times \sqrt{\sum_{k=1}^{m} \|D'_{nk} - D'_{mk}\|^2}$$
$$= o(\sigma_n^2).$$

Let $V_n(t) = \sigma_n^{-2} \sum_{k=1}^{m} E(D'^2_{nk}|\mathcal{F}_{k-1})$. Then $V_n(t) - t \xrightarrow{p} 0$. Let $I \geq 2$ be an integer. Observe that $\sup_{t \leq 1} |V_n(t) - t| \leq \max_{i \leq I} |V_n(i/I) - i/I| + 1/I$. By first letting $n \to \infty$ and then $I \to \infty$, (13) follows. □

THEOREM 2. *Relation* (10) *holds iff there is a martingale approximation for which* (11) *and* (12) *hold.*

PROOF. Suppose first that there is a martingale approximation (3) for which (11) and (12) hold. By Lemma 3, assume without loss of generality that the martingale approximation is defined by (6). Then, it suffices to establish (10) for all subsequences $n_r$, $r \geq 1$, that increase to $\infty$ sufficiently fast as $r \to \infty$. Observe that $D_{nk}$, $k = 1, 2, \ldots$, are martingale differences with respect to $P^x$ for a.e. $x(\pi)$ by the Markov property. If $n_r \to \infty$ sufficiently quickly as $r \to \infty$, then (12) and (13) both hold with convergence in probability replaced by convergence w.p.1 $(P)$, and $\lim_{n\to\infty}(S_n - M_{nn})/\sigma_n = 0$ w.p.1, too. So, these relations hold w.p.1 $(P^x)$ for a.e. $x(\pi)$. Then, for a.e. $x(\pi)$, $\lim_{r\to\infty} F_{n_r}(x; z) = \Phi(z)$ for all $z$, by the martingale central limit theorem applied conditionally given $X_0 = x$, and (10) holds (along the subsequence) by the bounded convergence theorem.

The converse will be deduced from Theorem 2 of Gänssler and Häeusler (1979), that provides necessary conditions for the functional version of the martingale central limit theorem. If (10) holds, then so does (4), by Lemma 2; and then there is a stationary martingale approximation, by Theorem 1. So, the issues are (11) and (12). Let $\mathbb{B}$ denote a standard Brownian motion. Then, since the process is stationary and $S_n^*$ is asymptotically normal given $X_0$,

$$\frac{1}{\sigma_n}[S_{\lfloor nt_1 \rfloor}, S_{\lfloor nt_2 \rfloor} - S_{\lfloor nt_1 \rfloor}, \ldots, S_{\lfloor nt_k \rfloor} - S_{\lfloor nt_{k-1} \rfloor}]$$
$$\Rightarrow \quad [\mathbb{B}_{t_1}, \mathbb{B}_{t_2} - \mathbb{B}_{t_1}, \ldots, \mathbb{B}_{t_k} - \mathbb{B}_{t_{k-1}}]$$

for every choice of $0 < t_1 < t_2 < \cdots < t_k \leq 1$, where $\Rightarrow$ denotes convergence in distribution. For example, if $k = 2$, $0 < s < t < 1$, and $m = \lfloor nt \rfloor - \lfloor ns \rfloor$, then

$$\left| P[S_{\lfloor ns \rfloor} \leq \sigma_n y, S_{\lfloor nt \rfloor} - S_{\lfloor ns \rfloor} \leq \sigma_n z] - \Phi\left(\frac{y}{\sqrt{s}}\right)\Phi\left(\frac{z}{\sqrt{t-s}}\right) \right|$$
$$\leq \int_{\mathcal{X}} \left| F_m\left(x; \frac{\sigma_n z}{\sigma_m}\right) - \Phi\left(\frac{z}{\sqrt{t-s}}\right) \right| \pi\{dx\}$$
$$+ \Phi\left(\frac{z}{\sqrt{t-s}}\right) \left| P[S_{\lfloor ns \rfloor} \leq \sigma_n y] - \Phi\left(\frac{y}{\sqrt{s}}\right) \right|,$$

which approaches zero as $n \to \infty$ since $\sigma_m/\sigma_n \to \sqrt{t-s}$. Next let

(14) $$\mathbb{M}_n(t) = \frac{1}{\sigma_n} \sum_{k \leq nt} D_{nk}$$

for $0 \leq t < 1$, and $\mathbb{M}_n(1) = \mathbb{M}_n(1-)$. Then the finite-dimensional distributions of $\mathbb{M}_n$ converge to those of $\mathbb{B}$, since $|S_{\lfloor nt \rfloor} - M_{\lfloor nt \rfloor}|/\sigma_n \xrightarrow{p} 0$ for each $0 < t < 1$; and since $E[\mathbb{M}_n(t)^2] \sim nt E(D_{n1}^2)/\sigma_n^2 \to t$, it follows that each $\mathbb{M}_n(t)^2$, $n \geq 1$, is uniformly integrable for each $0 < t \leq 1$. It then follows from the martingale inequality that $\mathbb{M}_n$ is tight in $D[0,1]$. So, $\mathbb{M}_n$ converges to $\mathbb{B}$ in $D[0,1]$; and relations (11) and (12) then follow from Theorem 2 of Gänssler and Häeusler (1979). □

EXAMPLE 1 (Continued). For linear processes, relations (11) and (12) follow from (4), which implies that $D_{nk} = \bar{b}_n \eta_k$ and that $|\bar{b}_n|$ is slowly varying, for the stationary martingale approximation constructed in the proof of Theorem 1. On the other hand, if $a_n = n^{-\beta}$, where $1/2 < \beta < 1$, then $S_n/\sigma_n$ is asymptotically standard normal, but (4) and (10) do not hold.

In the next corollary, let $\pi_1$ denote the joint distribution of $X_0$ and $X_1$, so that $\pi_1(B) = P[(X_0, X_1) \in B]$ for measurable $B \subseteq \mathcal{X}^2$; and let $H_n(x_0, x_1) = h_n(x_1) - Qh_n(x_0)$, so that $D_{nk} = H_n(X_{k-1}, X_k)$ in (6).

COROLLARY 2. *If* (4) *holds and* $H_n/\sqrt{\ell(n)} \to H \in L^2(\pi_1)$, *then* (10) *holds.*

PROOF. Let $D_{nk} = H_n(X_{k-1}, X_k)$ be the martingale approximation (6) and let $D'_{nk} = \sqrt{\ell(n)} H(X_{k-1}, X_k)$ and $M'_{nk} = D'_{n1} + \cdots + D'_{nk}$. Then the $D'_{nk}$ provide another stationary martingale approximation, since $\|M_{nn} - M'_{nn}\|^2 = n\|D_{n1} - D'_{n1}\|^2 = n\|H_n - \sqrt{\ell(n)} H\|^2 = o(\sigma_n^2)$. Moreover, the $D'_{nk}$ satisfy (11) and (12). For example,

$$\frac{1}{\sigma_n^2} \sum_{k \leq nt} E(D'^2_{nk}|X_k) = \frac{1}{n} \sum_{k \leq nt} E[H(X_{k-1}, X_k)^2 | X_{k-1}] \to tE[H(X_0, X_1)^2]$$

by the ergodic theorem; and $E[H(X_0, X_1)^2] = \|H\|^2 = 1$, since $\|H_n\|^2 \sim \ell(n)$, by Lemma 1. Condition (12) may be obtained similarly. □

To relate the condition in Corollary 2 to the sums $S_n$, first observe that $H_n/\sqrt{\ell(n)}$ converges in $L^2(\pi_1)$ iff $D_{n1}/\sqrt{\ell(n)}$ converges in $L^2(P)$ and next that $D_{n1}$ is the average of $E(S_k|X_1) - E(S_k|X_0)$ over $k = 1, \ldots, n$. It is not difficult to see that if $[E(S_n|X_1) - E(S_n|X_0)]/\sqrt{\ell(n)}$ converges in $L^2(P)$, then so does $D_{n1}/\sqrt{\ell(n)}$. Woodroofe (1992) shows how the condition of Corollary 2 can be related to the Fourier coefficients of $g$ when $X_k$ is a Bernoulli or Lebesgue shift process.

EXAMPLE 2. For a quasi-linear process (8), $D_{n1} = \bar{b}_n(X_0)\eta_1$. So, if $\ldots \eta_{-1}, \eta_0, \eta_1, \ldots$ are i.i.d., $\sigma_{n,1} = o(\sigma_{n,2})$, and $\bar{b}_n/\sqrt{\ell(n)} \to b \neq 0$ in $L^2(\pi)$, the $D_{n1}/\sqrt{\ell(n)}$ converges in $L^2(P)$ and, therefore, (11) and (12) both hold.

3.1. *Strong mixing processes.* Many classical results concerning asymptotic normality for stationary processes require strong mixing conditions; see, for example, Peligrad (1986, 1996). Here we show how the strong mixing assumption is related to our main condition (4). Let $X_n = (\ldots, \xi_{n-1}, \xi_n)$ and $S_n = \xi_1 + \cdots + \xi_n$, where $(\xi_i)_{i \in \mathbb{Z}}$ is a stationary sequence that is strong mixing; that is,

$$\alpha_n := \sup_{A \in \mathcal{F}_0, B \in \mathcal{G}_n} |P(A \cap B) - P(A)P(B)| \to 0$$

as $n \to \infty$, where $\mathcal{F}_n = \sigma(\ldots, \xi_{n-1}, \xi_n)$ and $\mathcal{G}_n = \sigma(\xi_n, \xi_{n+1}, \ldots)$.

LEMMA 4. *If $F$ and $G$ are two distribution functions and $\varepsilon > 0$, then there are continuous functions $w_1, \ldots, w_m$, depending only on $\varepsilon$ and $G$, for which $|w_i| \leq 1$ and $\int_{\mathbb{R}} w_i \, dG = 0$ for all $i$ and*

$$\Delta(G, F) \leq \varepsilon + \max_{i \leq m} \left| \int_{\mathbb{R}} w_i \, dF - \int_{\mathbb{R}} w_i \, dG \right|.$$

PROOF. The proof consists of first finding $a$ and $b$ for which $G(a) + 1 - G(b) \leq \varepsilon$, then partitioning $[a, b]$ into $a = x_0 < x_1 < \cdots < x_m = b$, where $x_i - x_{i-1} \leq \varepsilon/2$, constructing piecewise linear functions $u_i$ for which $u_i(x) = 1$ for $x \leq x_{i-1}$ and $u_i(x) = 0$ for $x \geq x_i$, and then letting $w_i = u_i - \int_R u_i \, dG$. The details are omitted. □

PROPOSITION 1. *Assume that $(\xi_n)_{n \in \mathbb{Z}}$ is a strong mixing process with mean 0 and finite variance. Then $S_n^* \Rightarrow \Phi$ implies* (10), *and consequently* (4).

PROOF. By Lemma 4, it suffices to show that

$$\int_{\mathcal{X}} \left| \int_{\mathbb{R}} w(z) F\{x; dz\} \right| \pi\{dx\} \to 0$$

as $n \to \infty$ for all continuous $w : \mathbb{R} \to [-1, 1]$ for which $\int_{\mathbb{R}} w \, d\Phi = 0$; and since the inner integral is just $E^x[w(S_n^*)]$, it suffices to show that $E|E[w(S_n^*)|X_0]| \to 0$ as $n \to \infty$ for all such $w$. To see this, let $m = m_n$ be a sequence for which $m \to \infty$ and $S_m/\sigma_n \xrightarrow{p} 0$; and let $\tilde{S}_n = (S_{n+m} - S_m)/\sigma_n$. Further, let $w : \mathbb{R} \to [-1, 1]$ be a continuous function for which $\int_{\mathbb{R}} w \, d\Phi = 0$ and let $w_n(x) = E^x[w(S_n^*)]$ and $\tilde{w}_n(x) = E^x[w(\tilde{S}_n)]$. Then $E[w_n(X_0)] = E[w(S_n^*)] \to 0$, since $S_n^* \Rightarrow \Phi$; $E|w_n(X_0) - \tilde{w}_n(X_0)| \leq E|w(S_n^*) - w(\tilde{S}_n)| \to 0$, since $\tilde{S}_n - S_n^* \xrightarrow{p} 0$ as $n \to \infty$; and

$$E|\tilde{w}_n(X_0)|^2 = \int \tilde{w}_n(X_0) w(\tilde{S}_n) \, dP \leq E[w(\tilde{S}_n)]^2 + 4\alpha_m \to 0,$$

by standard mixing inequalities [see, e.g., Hall and Heyde (1980), page 277]. So, $E|w_n(X_0)| \to 0$ as $n \to \infty$ as required. □

**4. An invariance principle.** Let

$$\mathbb{B}_n(t) = \frac{1}{\sigma_n} S_{\lfloor nt \rfloor}$$

for $0 \leq t < 1$, $\mathbb{B}_n(1) = \mathbb{B}_n(1-)$, where $\lfloor x \rfloor$ denotes the greatest integer that is less than or equal to $x$. If (10) holds, then the finite-dimensional distributions of $\mathbb{B}_n$ converge to those of standard Brownian motion $\mathbb{B}$, and $\mathbb{M}_n$ converges in distribution to $\mathbb{B}$ in the space $D[0, 1]$, both from the proof of Theorem 2. Relations (4) and (10) do not imply that $\mathbb{B}_n$ converges in distribution to $\mathbb{B}$ in $D[0, 1]$, however.

EXAMPLE 3. Let $G$ be a symmetric distribution function for which

$$1 - G(y) \sim \frac{1}{y^2 \log^{3/2}(y)}$$

as $y \to \infty$. Let $\ldots, \eta_{-1}, \eta_0, \eta_1, \ldots \sim \Phi$ and $\ldots, Y_{-1}, Y_0, Y_1, \ldots \sim G$ be independent random variables. Let $a_0 = 0, a_1 = 1/\log(2)$ and $a_k = 1/\log(k+1) - 1/\log(k)$ for $k \geq 2$, as in Example 1. Define $\xi_k$ by (8); let $\xi_k' = \xi_k + Y_k - Y_{k-1}$; and let $S_n = \xi_1 + \cdots + \xi_n$ and $S_n' = \xi_1' + \cdots + \xi_n'$. Then (4), (11) and (12) hold for both $S_n$ and $S_n'$ with $\sigma_n^2 \sim n/\log^2(n)$. In this example,

$$\frac{1}{\sigma_n} \max_{k \leq 1} |Y_k - Y_0| \to \infty$$

in probability, so that $\mathbb{B}_n$ and $\mathbb{B}_n'$ cannot both converge to $\mathbb{B}$.

In Theorem 3 and Corollary 3, we consider the special case in which $\sigma^2 = \lim_{n \to \infty} \sigma_n^2/n$ exists. These results improve Theorem 2 and Corollary 4 in Maxwell and Woodroofe (2000) by imposing a weaker condition as well as by obtaining a stronger result. The heart of the matter is whether there is a martingale approximation for which $\max_{k \leq n} |S_k - M_{nk}|/\sqrt{n} \to 0$ in probability. This question is addressed first. Two lemmas are needed.

LEMMA 5. *Suppose that, for some $q > 1$,*

(15) $$\|E(S_n|X_0)\| = o(\sqrt{n} \log^{-q} n).$$

*Then there is a martingale $M_1, M_2, \ldots$ with stationary increments for which $\|S_n - M_n\| = o(\sqrt{n} \log^{1-q} n)$.*

PROOF. Recall the construction $D_{nk}$ and $M_{nk} = D_{n1} + \cdots + D_{nk}$ from (6) and also that $\max_{k \leq n} \|S_k - M_{nk}\| \leq 3 \max_{k \leq n} \|E(S_k|X_0)\|$. Thus, $\max_{k \leq n} \|S_k - M_{nk}\| = o[\sqrt{n} \log^{-q}(n)]$ in the present context. So, if $m \geq 2$ and $m \leq n \leq 2m$, then $\|M_{nm} - M_{mm}\| = o[\sqrt{m} \log^{-q}(m)]$. Since $\|M_{nm} - M_{mm}\|^2 = m\|D_{n1} - D_{m1}\|^2 = m\|H_n - H_m\|^2$, it then follows that

(16) $$\sum_{k=j}^{\infty} \|H_{2^k} - H_{2^{k-1}}\| \leq \sum_{k=j}^{\infty} o[\log^{-q}(2^k)] = o[\log^{1-q}(2^j)].$$

It follows that $H_{2^k}$ has a limit $H$, say, in $L^2(\pi_1)$ and that $\|H - H_m\| = o[\log^{1-q}(m)]$. Letting $D_k = H(X_{k-1}, X_k)$ and $M_n = D_1 + \cdots + D_n$, the lemma then follows from $\|S_n - M_n\| \leq \|S_n - M_{nn}\| + \sqrt{n}\|H_n - H\|$. □

LEMMA 6. *Let $Y_k$, $k \in \mathbb{Z}$, be a second-order stationary process with mean 0 and let $T_n = Y_1 + \cdots + Y_n$. Then*

$$E\left[\max_{k \leq n} T_j^2\right] \leq d \sum_{j=0}^{d} 2^{d-j} \|T_{2^j}\|^2,$$

where $d = \lceil \log_2(n) \rceil$, the least integer that is greater than or equal to $\log_2(n)$.

PROOF. The proof uses a simple chaining argument and appears in Doob [(1953), page 156] for uncorrelated random variables. Briefly, any integer $k \leq n$ may be written as $k = 2^{r_1} + \cdots + 2^{r_j}$, where $0 \leq r_j < \cdots < r_1 \leq d$. So,

$$|T_k|^2 = \left| \sum_{i=1}^{j} (T_{2^{r_1}+\cdots+2^{r_i}} - T_{2^{r_1}+\cdots+2^{r_{i-1}}}) \right|^2$$

$$\leq j \sum_{i=1}^{j} |T_{2^{r_1}+\cdots+2^{r_i}} - T_{2^{r_1}+\cdots+2^{r_{i-1}}}|^2,$$

where an empty sum is to be interpreted as 0, and

$$\max_{k \leq n} |T_k|^2 \leq d \sum_{j=0}^{d} \sum_{i=1}^{2^{d-j}} |T_{i2^j} - T_{(i-1)2^j}|^2,$$

from which the lemma follows by stationarity. □

THEOREM 3. *Let $R_n = S_n - M_n$, where $M_n$ is as in Lemma 5. If $g \in L^p$ for some $p > 2$ and (15) holds for $q \geq 2$, then $\sigma^2 = \lim_{n \to \infty} \sigma_n^2/n$ exists, and*

$$(17) \qquad \lim_{n \to \infty} P\left[ \max_{j \leq n} |R_j| \geq \varepsilon \sqrt{n} \right] = 0$$

*for each $\varepsilon > 0$; and if (15) holds for some $q > 5/2$, then*

$$\lim_{n \to \infty} P^x\left[ \max_{j \leq n} |R_j| \geq \varepsilon \sqrt{n} \right] = 0$$

*for a.e. $x(\pi)$ for each $\varepsilon > 0$.*

PROOF. Let $\gamma = 1/4 - 1/(2p) > 0$, where $p$ is as in the statement of the theorem, $a = a_m = \lceil 2^{m\gamma} \rceil$, and $b = b_m = \lceil 2^{m(1-\gamma)} \rceil$. Then

$$\max_{j \leq 2^m} |R_j| \leq \max_{1 \leq k \leq b} \left[ |R_{ak}| + \max_{0 \leq j \leq a} |R_{ak+j} - R_{ak}| \right].$$

Here,

$$\max_{0 \leq j \leq a} |R_{ak+j} - R_{ak}|$$

$$\leq \max_{0 \leq j \leq a} |M_{ak+j} - M_{ak}| + \max_{0 \leq j \leq a} |S_{ak+j} - S_{ak}|$$

$$\leq \max_{0 \leq j \leq a} |M_{ak+j} - M_{ak}| + a \max_{j \leq 2^m} |g(X_j)|$$

for each $k$. So,

(18)
$$P^x\left[\max_{j\leq 2^m}|R_j|\geq 3\varepsilon\sqrt{2^m}\right]$$
$$\leq P^x\left[\max{}^*\frac{|M_k-M_j|}{\sqrt{2^m}}\geq\varepsilon\right]$$
$$+P^x\left[\max_{j\leq 2^m}\frac{|g(X_j)|}{\sqrt{2^m}}\geq\frac{\varepsilon}{a}\right]+P^x\left[\max_{k\leq b}\frac{|R_{ak}|}{\sqrt{2^m}}\geq\varepsilon\right],$$

where $\max^*$ runs over all pairs $(j,k)$ such that $1\leq j, k\leq 2^m$ and $|k-j|\leq a$. The first term clearly tends to 0 for a.e. $x(\pi)$, by the functional martingale central limit theorem. The second term in (18) also converges to 0 for a.e. $x(\pi)$ by the Borel–Cantelli lemma, since

$$\int_{\mathcal{X}}P^x\left[\max_{j\leq 2^m}\frac{|g(X_j)|}{\sqrt{2^m}}\geq\frac{\varepsilon}{a}\right]\pi\{dx\}\leq\frac{a^p}{\varepsilon^p}2^{m(1-p)}E|g(X_1)|^p,$$

and the right-hand side is summable over $m$ (recalling that $a=\lceil 2^{\gamma m}\rceil$ and observing that $p\gamma+1-p<0$). Similarly, for the third term on the right-hand side of (18),

$$\int_{\mathcal{X}}P^x\left[\max_{k\leq b}\frac{|R_{ak}|}{\sqrt{2^m}}\geq\varepsilon\right]\pi\{dx\}=P\left[\max_{k\leq b}\frac{|R_{ak}|}{\sqrt{2^m}}\geq\varepsilon\right]$$
$$\leq\frac{1}{\varepsilon^2}E\left[\max_{k\leq b}\frac{|R_{ak}|}{\sqrt{2^m}}\right]^2,$$

and, letting $d=\lceil\log_2(b)\rceil$,

$$E\left[\max_{k\leq b}\frac{|R_{ak}|}{\sqrt{2^m}}\right]^2\leq\frac{d}{2^m}\sum_{i=0}^{d}2^{d-i}\|R_{a2^i}\|^2$$
$$\leq\frac{d}{2^m}\sum_{i=0}^{d}2^{d-i}\frac{o(a2^i)}{\log^{2(q-1)}(a2^i)}$$
$$=\frac{abd}{2^m}o\left[\frac{1}{m^{2q-3}}\right]=o(m^{4-2q}),$$

by Lemmas 5 and 6. Relation (17) follows immediately, since $ab=O(2^m)$ and $d=O(m)$; and if $q>5/2$, then $o(m^{4-2q})$ is summable and $P^x[\max_{k\leq b}|R_{ak}|\geq\varepsilon\sqrt{2^m}]\to 0$ for a.e. $x$, by the Borel–Cantelli lemma. $\square$

Now let $G_n$ and $\Psi$ be the distributions of $\mathbb{B}_n$ and Brownian motion in $D[0,1]$, and let $\Delta$ denote the Prokhorov metric for $D[0,1]$.

COROLLARY 3. *If* (15) *holds for some* $q \geq 2$ *and* $0 < \sigma^2 < \infty$, *then*

$$\lim_{n \to \infty} \int_{\mathcal{X}} \Delta[\Psi, G_n(x; \cdot)] \pi\{dx\} = 0;$$

*and if* $q > 5/2$ *in* (15), *then* $\lim_{n \to \infty} \Delta[\Psi, G_n(x; \cdot)] = 0$ *for a.e.* $x(\pi)$.

PROOF. Let $K_n(x; \cdot)$ be the distribution of $\mathbb{M}_n$ in $D[0,1]$. Then $K_n(x; \cdot) \Rightarrow \Psi$ as $n \to \infty$ for a.e. $x(\pi)$, by the functional central limit theorem, and

$$\Delta[\Psi, G_n(x; \cdot)] \leq \Delta[\Psi, K_n(x; \cdot)] + P^x\left[\max_{k \leq n} |R_k| \geq \varepsilon \sigma_n\right] + \varepsilon$$

for each $\varepsilon > 0$. The case $q > 5/2$ follows immediately, and the case $2 \leq q \leq 5/2$ from $\int_{\mathcal{X}} P^x[\max_{k \leq n} |R_k| \geq \varepsilon \sigma_n] \pi\{dx\} = P[\max_{k \leq n} |R_k| \geq \varepsilon \sigma_n]$. $\square$

COROLLARY 4. *If* (15) *holds for some* $q \geq 2$ *and* $\sigma^2 = 0$, *then* $\max_{k \leq n} |S_k|/\sqrt{n} \xrightarrow{p} 0$; *and if* $q > 5/2$, *then* $\lim_{n \to \infty} P^x[\max_{k \leq n} |S_k| \geq \varepsilon \sqrt{n}] = 0$ *for a.e.* $x(\pi)$ *for each* $\varepsilon > 0$.

PROOF. In this case $S_k = R_k$. $\square$

REMARK 1. A simple sufficient condition for (15) is

(19) $$\|E[g(X_n)|X_0]\| = \mathcal{O}(n^{-1/2} \log^{-q} n).$$

However, (15) allows processes of the form (8) with $a_n = n^{-\beta}(-1)^n$ for $n \geq 1$, where $1/2 < \beta < 1$. In this case (19) is violated. Wu (2002) derived central limit theorems for processes of this sort whose covariances are summable but not absolutely summable. A typical example is the Gegenbauer process which exhibits long-range dependence and has oscillatory covariances [Beran (1994)].

**Acknowledgment.** The authors gratefully acknowledge the helpful and constructive comments provided by the anonymous referee.

Department of Statistics
University of Chicago
5734 S. University Avenue
Chicago, Illinois
USA
e-mail: wbwu@galton.uchicago.edu

Department of Statistics
University of Michigan
462 West Hall
Ann Arbor, Michigan
USA
e-mail: michaelw@umich.edu